\documentclass[11pt, a4paper]{article} 
\usepackage{amsfonts,fancyhdr} 
\def\summe#1#2{
	\mathop{{\sum}^{\ne}}
	\limits_{#1\hspace*{0.6em}}^{#2\hspace*{0.6em}}}

\def\hut#1{\widehat{#1}}

\newcommand{\E}{\mathbf{E}}
\newcommand{\eins}{\mathbf{1}}
\newcommand{\Var}{\mathbf{V}}
\newcommand{\Wk}{\mathbf{P}}
\newcommand{\R}{\mathbb{R}}

\begin{document}

\renewcommand{\title}[1]{
{\LARGE\bf\begin{center}#1\end{center}}
\vspace{0.5cm}
{\small
\begin{flushleft}
	\today\\
	\vspace{0.2cm}
	Martin Snethlage\\
	Graduiertenkolleg ``R\"aumliche Statistik''\\
	TU Bergakademie Freiberg\\
	Bernhard-von-Cotta-Str. 2\\
	D-09596 Freiberg\\
	\vspace{0.1cm}
	snethlag@grad.tu-freiberg.de
\end{flushleft}
\vspace{0.1cm}
}}

\title{Is Bootstrap Really Helpful in Point Process Statistics?}

\begin{abstract}
There are some papers which 
describe the use of bootstrap techniques in point process
statistics. The aim of the present paper is to show that the form in which 
bootstrap is used there is dubious. In case of variance estimation of 
pair correlation function estimators
the used bootstrap techniques lead to results
which can be obtained simpler without simulation;
furthermore, they differ from the desired results.
The problem to obtain
confidence regions for the intensity function
of inhomogeneous Poisson processes can be easily solved without
bootstrap techniques.
\end{abstract}


\section{Introduction}
\label{intro}
Recently, bootstrap is a popular tool in many
branches of statistics, also
for stochastic processes. Thus it is natural to ask 
whether bootstrap techniques
could be helpful also in point process statistics.
Indeed, some authors have developed statistical procedures
using bootstrap techniques, 
see e.g. \cite{Barrow84}, \cite{Borner92}, \cite{Cowling96}
and \cite{Ling86}.
All these papers deal with the estimation of the accuracy
of estimators of point process characteristics.
In the first three papers
estimation of variance of
pair correlation function estimators is treated.
The last one presents a procedure to determine confidence
regions for the intensity function of an
inhomogeneous Poisson process.

The fundamental idea of bootstrap 
to resample given data to obtain `new' pseudo data
appears also in statistics of stochastic processes,
in particular in the analysis of time series, see e.g. \cite{Kuensch89}. 
In some variants of the method, called the blockwise
bootstrap, the time series is partitioned
into several parts, which are then resampled. 
A similar idea is also applied in \cite{Hall85}
in the statistical analysis of a planar random set.
Clearly, the partioning procedure
can also be adapted to point process statistics. 
However, partition
can destroy point structures or add new
artificial structures to the point pattern.
In the one-dimensional case the error resulting from this
loss of information
may be still acceptable, but in higher dimensions
it will be serious.
Thus in spatial point process statistics another
method is used which is quite 
similar to the application of bootstrap
in case of classical statistics:
the points of the process (including 
their places, which are assumed to be pairwise different)
are resampled.
The pseudo
pattern then consists of $n$ points $x^*_1,\ldots,x^*_n$
which are obtained by sampling randomly with replacement $n$ times from 
the original data $\{x_1,\ldots,x_n\}$.

Naturally, the pseudo patterns generated by
this method
have always multiple points.
Thus they have a character different 
to that of the original, which does not have multiple points.

Consequently, it would be surprising if 
quantities of such point processes 
would produce good estimators for quantities of the original point process.

This paper analyses the pointwise resampling technique for some examples
of point process statistics.
Section \ref{var} discusses the main ideas of 
the paper
\cite{Borner92} which presents
a procedure for estimating the standard error
of an estimator of the pair correlation function.
In Section \ref{conf} a method (drawn from
\cite{Cowling96}) 
to determine confidence intervals
for the intensity function of an inhomogeneous
Poisson process is considered.
Finally, an easier method is presented
which yields confidence regions for
the intensity function of an inhomogeneous
Poisson process without bootstrap.

\section{Variance of estimators of the product density}
\label{var}
This part discusses the main ideas of  \cite{Borner92}, 
where bootstrap techniques are used to
approximate the standard error of a pair correlation function estimator. 
The calculations are presented in an abridged form;
the complete calculations are given in the Appendix.

\subsection{Fundamentals}

Let $\Phi$ be a stationary and isotropic 
point process, see, for example, \cite{Stoyan95}
for definitions.
A standard second order characteristic of $\Phi$ 
is the product density function $\varrho^{(2)}(r)$.
This function can be interpreted heuristically
as follows.
If $B_1$ and $B_2$ are two infinitesimally small disjoint
Borel sets of volumes d$V_1$ and d$V_2$ and if $x_1\in B_1$ and
$x_2\in B_2$ are points of distance $\Vert x_1-x_2\Vert=r$ then
$\varrho^{(2)}(r)\textrm{d}V_1\textrm{d}V_2$
is the probability that $\Phi$ has a 
point in each of $B_1$ and $B_2$. 
A simple estimator of $\varrho^{(2)}$
without any border correction is given by
$$
	\hut{\varrho}(r)=\frac{1}{2\pi r\nu(W)}\summe{x,y\in\Phi\cap W}{}
	K(r-\Vert x-y\Vert).
$$
The summation goes over all point pairs
with different members,
$W$ denotes the window of observation and
$K$ is a kernel function.

This situation can be generalized to the case of any
`two-point estimator'
\begin{equation}\label{pop2}
	\hut{\theta}=\summe{x,y\in\Phi}{}f(x,y)
\end{equation}
with $f$ being symmetrical 
in its arguments and 
of the form 
$$
	f(x,y)=\eins_W(x)\eins_W(y)h(x,y)
$$
with some function $h$. 
As the special form of $f$ leading to $\hut{\varrho}(r)$
is unimportant, the
following calculations are carried out for a general $\hut{\theta}$.

The quantity of interest is the variance of $\hut{\theta}$
which is given by
\begin{equation}\label{vartheta}
	\Var\hut{\theta}=\E\hut{\theta}^2-(\E\hut{\theta})^2
	= s_4+4s_3+2s_2-(\E\hut{\theta})^2
\end{equation}
with
$$
	s_i=\int\varrho^{(i)}(x_1,\ldots,x_i)
	f(x_1,x_2)f(x_{i-1},x_i)\,\textrm{d}x_1\ldots\textrm{d}x_i,
$$
where $\varrho^{(i)}$ is the $i$th order product density function
of $\Phi$, see the Appendix.

\subsection{Bootstrap version of $\hut{\theta}$}

Assume that a sample of $\Phi$ is given which consists of $n$ points 
$x_1,\ldots,x_n$ in
the observation window $W$.
It is resampled $N$ times to obtain $N$ `new' point patterns. 
Each pseudo pattern consists of $n$ points $x^*_1,\ldots,x^*_n$
which are obtained by sampling randomly with replacement $n$ times from 
$\{x_1,\ldots,x_n\}$.
Thus it happens that in the pseudo samples some points 
of the original point pattern do not occur
while others
occur twice or even more.
Let the number of occurrences of 
$x_i$ in the $k$th sample be $w_k(i)$.
Then the $k$th sample can be represented by
the vector $w_k=(w_k(1),\ldots,w_k(n))$
which has a multinomial
distribution.
This distribution depends only on $n$. 
In the limiting case $n\to\infty$ the components $w_k(i)$ of $w_k$
are independent and Poisson distributed with mean $\mu=1$.

The bootstrap estimate for the $k$th pseudo sample is
$$
	\hut{\theta^*_k}=\summe{i,j=1}{n}
	f(x_i,x_j)w_k(i)w_k(j),\quad k=1,\ldots,N
$$
where the summation goes over all pairs $(i,j)$ with $i\neq j$. 
The variance of $\hut{\theta}$ is estimated by the usual 
variance estimator corresponding to the $\hut{\theta^*_k}$'s,
$$
	\hut{v^*_N}=\frac{1}{N-1}\sum_{k=1}^{N}
	\left(\hut{\theta^*_k}-\frac{1}{N}
	\sum_{i=1}^{N}\hut{\theta^*_i}\right)^2 .
$$

Since the $\hut{\theta^*_k}$ are (conditionally on
$x_1,\ldots,x_n$) independent and identically distributed,
it is 
\begin{eqnarray}
	\lim_{N\to\infty}\hut{v^*_N} & 
	= & \mathbf{V}\hut{\theta^*_1}\nonumber\\
	& = & \E\hut{\theta^*_1}^2-
	(\E\hut{\theta^*_1})^2\nonumber\\
	& = & \alpha_4\summe{i,j,k,l=1}{n}
	f(x_i,x_j)f(x_k,x_l)\label{vartheta*}\nonumber\\
	& & {}+4\alpha_3\summe{i,j,k=1}{n}
	f(x_i,x_j)f(x_i,x_k)\\
	& & {}+2\alpha_2\summe{i,j=1}{n}
	\left(f(x_i,x_j)\right)^2\nonumber
\end{eqnarray}
with
\begin{eqnarray*}
	\alpha_4 & = & \left[\E w_1(1)w_1(2)w_1(3)w_1(4)-
	\left(\E w_1(1)w_1(2)\right)^2\right]\\
	\alpha_3 & = & \left[\E(w_1(1))^2 w_1(2)w_1(3)-
	\left(\E w_1(1)w_1(2)\right)^2\right]\\
	\alpha_2 & = & \left[\E(w_1(1)w_1(2))^2-
	\left(\E w_1(1)w_1(2)\right)^2\right],\nonumber\\
\end{eqnarray*}
where the expectations are conditionally on fixed
$x_1,\ldots,x_n$.
All the $\alpha_i$ can be calculated numerically and depend only on $n$
(see the Appendix).
Thus the result of the whole bootstrap procedure for $N\to\infty$
can be simply obtained by direct computation.

\subsection{Expectation of $\hut{v^*_N}$}\label{sack}
The futility
of $\hut{v^*_N}$ is demonstrated
by the fact that it
does neither estimate what is hoped (the variance of $\hut{\theta}$)
nor a multiple with a fixed factor. 
To show this, the unconditional expectation of $\hut{v^*_N}$
is determined, see the Appendix. Since the
result is not very transparent, here an 
approximation is given which makes it possible to characterize
the quality of $\hut{v^*_N}$. 

Assume that
the $w_k(i)$ are independent and 
Poisson distributed with parameter $\mu=1$;
this simplifying assumption is 
exact in the limiting case $n\to\infty$, see above.
This leads to a result which is close
to the exact value for large $n$ and is easy to interpret. 
By the way, the simplification is equivalent
to replacement of $n$ by $n^*$ in each pseudo
sample where $n^*$ is a 
Poisson distributed number with mean $\mu=n$. In this scheme  
each pseudo
sample consists of a random number of points.

The result is
\begin{equation}\label{result}
	\lim_{N\to\infty}\E\hut{v^*_N} = 
	\E\lim_{N\to\infty}\hut{v^*_N} = 
	4s_3+6s_2,
\end{equation}
see the Appendix, while the desired result, given by (\ref{vartheta}),
is
$$
\Var\hut{\theta}=s_4+4s_3+2s_2-\left(\E\hut{\theta}\right)^2.
$$

{\em Remark:} The formulae suggest that the bootstrap result
(\ref{result}) can considerably differ from the true
variance of $\hut{\theta}$. Nevertheless, the bootstrap
procedure may make sense. 
In some cases $s_4$ converges to 
$\left(\E\hut{\theta}\right)^2$ with growing
$W$ and $s_3$ is small compared with
$s_2$. Then the
bootstrap result (\ref{result}) may
approximate three times the true variance,
see \cite{Borner92}.

\section{Confidence regions for the intensity function 
	of an inhomogeneous Poisson process}
\label{conf}
The paper \cite{Cowling96} presents a procedure
which uses bootstrap techniques
to determine confidence regions for
the intensity function of an
inhomogeneous Poisson process. 
The confidence regions are estimated using a kernel estimator.
The following discusses the main idea of that
paper and shows that, as above, it is not necessary to carry out
the bootstrap procedure.

\subsection{Fundamentals}

For simplicity, the following calculations 
are carried out for an one-dimensional
point process, but
they could be easily generalized to 
higher-dimensional processes.

Consider an inhomogeneous Poisson point
process $\Phi$ with unknown intensity function 
$\lambda(x)$ in the interval $(0,1)$,
with points $0<x_1\leq x_2\leq\ldots\leq x_n<1$.
A kernel estimator for
$\lambda(x)$ is used as
$$
	\hut{\lambda}(x)=
	\frac{1}{h}\sum_{i=1}^nK\left(\frac{x-x_i}{h}\right),\quad x\in(0,1),
$$
where $K$ is a kernel function and $h$ bandwidth
(see, for example, \cite{Diggle83}).

Define
$$
	T(x)=\frac{\hut{\lambda}(x)-\E\hut{\lambda}(x)}
	{\sqrt{\hut{\lambda}(x)}},\quad 0<x<1,
$$
and, for a given $\alpha$ with $0\leq\alpha\leq1$,
$$
	t_\alpha(x)=\min_{t\in\R^+}\left\{t:\Wk\left\{\left\vert 
	T(x)\right\vert\leq t\right\}\geq1-\alpha\right\},\quad 0<x<1.
$$
Then an estimate of a confidence region for $\lambda(x)$ 
of level $1-\alpha$ is the interval
$$
	C(x)=\left[\hut{\lambda}(x)-t_\alpha(x)\sqrt{\hut{\lambda}(x)},
	\hut{\lambda}(x)+t_\alpha(x)\sqrt{\hut{\lambda}(x)}\right],\quad 0<x<1,
$$
where the left border is set on 0 if it is negative.

\subsection{Bootstrap versions}

Since the distribution of $T$ is not available 
(because the intensity function is unknown) it is 
approximated by simulation of pseudo data, see \cite{Cowling96}.
A set of pseudo data is obtained 
by drawing $x^*_1,\ldots,x^*_{n^*}$
by sampling randomly
with replacement $n^*$ times from $\{x_1,\ldots,x_n\}$, where
$n^*$ has a Poisson distribution with mean $n$
(this is method 2 in \cite{Cowling96} and 
similar to the simplified case in Section \ref{sack}).
The number of occurrences of 
$x_i$ in the $k$th sample is a random variable, denoted 
as above by 
$w_k(i)$.
All the $w_k(i)$
are independent and Poisson distributed with mean $\lambda=1$ for
$i=1,\ldots,n$.

For given $\alpha$ with $0\leq\alpha\leq1$ the bootstrap versions 
of the quantities defined above are
\begin{eqnarray*}
	\hut{\lambda^*_k}(x) & = &
	\frac{1}{h}\sum_{i=1}^nK\left(\frac{x-x_i}{h}\right)w_k(i),\\
	T^*_k(x) & = & \frac{\hut{\lambda^*_k}(x)-\hut{\lambda}(x)}
	{\sqrt{\hut{\lambda^*_k}(x)}},\quad x\in(0,1),\\
	t_\alpha^*(x) & = & \min_{t\in\R^+}\left\{t:\Wk^*\left\{\left\vert 
	T^*(x)\right\vert\leq t\right\}\geq1-\alpha\right\},\\
	C^*(x) & = & \left[\hut{\lambda}(x)-t_\alpha^*(x)\sqrt{\hut{\lambda}(x)},
	\hut{\lambda}(x)+t_\alpha^*(x)\sqrt{\hut{\lambda}(x)}\right],
\end{eqnarray*}
where $\Wk^*(\,\cdot\,)=\Wk(\,\cdot\,\vert\{x_1,\ldots,x_n\})$ 
is the distribution conditionally on 
\linebreak$\{x_1,\ldots,x_n\}$.

The determination of $t_\alpha^*(x)$ can be carried out by simulation.
However, a faster and simpler possibility uses the well-known fact that the 
sum of independent
Poisson distributed random variables
is also Poisson distributed.
It is demonstrated here for the simple
rectangular kernel function
$$
K(x)=\frac{1}{2}\cdot\eins_{\lbrack-1,1\rbrack}(x).
$$
For other kernels, similar calculations are possible.
Let $p(x)$ be the number of observed points in the interval
$\left[x-h,x+h\right]$.
Then its bootstrap version $p^*(x)$ is a random variable
which is Poisson distributed with mean $p(x)$.
Its cumulative disribution function is denoted by $\mathbf{F}^*$.
Thus, for given $\alpha$,
\begin{eqnarray}
	t_\alpha^*(x) & = & \min_{t\in\R^+}\left\{t:\Wk^*\left\{\left\vert 
	T^*_1(x)\right\vert\leq t\right\}\geq1-\alpha\right\}\nonumber\\
	 & = & \min_{t\in\R^+}
	\bigg\{t:\Wk^*\bigg\{
	\Big\vert\sum_{i=1}^{p(x)}w_1(i)-a_{x}(t)
	\Big\vert\leq b_{x}(t)
	\bigg\}\geq1-\alpha\bigg\}\label{t*1}\label{popel4}\\
	 & = & \min_{t\in\R^+}\left\{t:
	\mathbf{F}^*\left(a_{x}(t)+b_{x}(t)\right)-
	\mathbf{F}^*\left(a_{x}(t)-b_{x}(t)\right)
	\geq1-\alpha\right\}\nonumber,
\end{eqnarray}
with 
\begin{eqnarray*}
	a_{x}(t) & = & p(x)+ht^2,\\
	b_{x}(t) & = & t\sqrt{2hp(x)+h^2t^2}
\end{eqnarray*}
The calculations are of an elementary nature
using that $\hut{\lambda_k^*}$ is equal to
$\frac{1}{2h}\sum^{p(x)}_{1=1}w_k(i)$.

\subsection{Confidence regions without bootstrap}

Of course, in case of an inhomogeneous
Poisson process it is easy
to build confidence regions without the bootstrap methodology.
Assume that the intensity function $\lambda(x)$ is
approximately linear in the interval $\left[x-h,x+h\right]$. Then
$2h\hut{\lambda}(x)$ using the rectangular kernel
is Poisson distributed with mean $2h\lambda(x)$.
(The rectangular kernel could be replaced
by another kernel; then the corresponding 
calculations become a bit more difficult.)
Therefore, known confidence regions for the Poisson
parameter can be used, see for example \cite{Sachs84}.
Thus it is easy to build the
desired confidence region for $\lambda(x)$.

This result corresponds to a general 
observation. If a parametric statistic 
problem is given, then
parametric estimators lead usually to better results than
bootstrap techniques.

\vspace{0.7cm}
{\small
\noindent{\em Acknowledgements:} 
\label{ackn}
I am most grateful to Dietrich Stoyan for his great encouragement
and for helpful discussions.
}

\section*{Appendix}

Here the derivation of some equations of Section \ref{var}
is given.

\subsection*{Equation (\ref{vartheta})}
For
$$
	\hut{\theta}=\summe{x,y\in\Phi}{}f(x,y)
	=\summe{x,y\in\Phi}{}\eins_W(x)\eins_W(y)h(x,y)
$$
it is
\begin{eqnarray*}
	\E\hut{\theta}^2 & = & \E\bigg(\summe{x,y\in\Phi}{}
	f(x,y)\bigg)^2\nonumber\\
	& = & \hspace{0.45cm}\E\summe{w,x,y,z\in\Phi}{}
	f(w,x)f(y,z)\nonumber\\
	& & {}+4\E\summe{x,y,z\in\Phi}{}f(x,y)f(x,z)\\
	& & {}+2\E\summe{x,y\in\Phi}{}
	\left(f(x,y)\right)^2\nonumber\\
	& = & \hspace{0.73cm}\int\varrho^{(4)}(x_1,x_2,x_3,x_4)f(x_1,x_2)f(x_3,x_4)
	\textrm{d}\, x_1\textrm{d}\, x_2\textrm{d}\, x_3\textrm{d}\, x_4\nonumber\\
	& & {}+4\int\varrho^{(3)}(x_1,x_2,x_3)
	f(x_1,x_2)f(x_1,x_3)\textrm{d}\, x_1\textrm{d}\, x_2\textrm{d}\, x_3\\
	& & {}+2\int\varrho^{(2)}(x_1,x_2)
	\left(f(x_1,x_2)\right)^2\textrm{d}\, x_1\textrm{d}\, x_2\nonumber\\
	& = & s_4+4s_3+2s_2\nonumber
\end{eqnarray*}
with
$$
	s_i=\int\varrho^{(i)}(x_1,\ldots,x_i)
	f(x_1,x_2)f(x_{i-1},x_i)\textrm{d}\, x_1\ldots\textrm{d}\, x_i.
$$

\subsection*{Equation (\ref{vartheta*})}
For
$$
	\hut{\theta^*_1}=\summe{i,j=1}{n}f(x_i,x_j)w_1(i)w_1(j)
$$
it is
\begin{eqnarray*}
	\E\hut{\theta^*_1}^2 & = & 
	\E\hut{\theta^*_1}^2-
	\E\hut{\theta^*_1}\hut{\theta^*_2}\nonumber\\
	& = & \hspace{0.5cm}\E
	\summe{i,j,k,l=1}{n}f(x_i,x_j)f(x_k,x_l)
	w_1(i)w_1(j)w_1(k)w_1(l)\nonumber\\
	& & {}+4\E\summe{i,j,k=1}{n}f(x_i,x_j)f(x_i,x_k)
	(w_1(i))^2w_1(j)w_1(k)\\
	& & {}+2\E\summe{i,j=1}{n}\left(f(x_i,x_j)\right)^2
	(w_1(i)w_1(j))^2\nonumber\\
	& = & \hspace{0.55cm}\summe{i,j,k,l=1}{n}f(x_i,x_j)f(x_k,x_l)
	\E w_1(i)w_1(j)w_1(k)w_1(l)\nonumber\\
	& & {}+4\summe{i,j,k=1}{n}f(x_i,x_j)f(x_i,x_k)
	\E(w_1(i))^2w_1(j)w_1(k)\\
	& & {}+2\summe{i,j=1}{n}\left(f(x_i,x_j)\right)^2
	\E(w_1(i)w_1(j))^2\nonumber\\
	& = & \hspace{0.5cm}\E w_1(1)w_1(2)w_1(3)w_1(4)
	\summe{i,j,k,l=1}{n}f(x_i,x_j)f(x_k,x_l)\nonumber\\
	& & {}+4\E(w_1(1))^2 w_1(2)w_1(3)
	\summe{i,j,k=1}{n}f(x_i,x_j)f(x_i,x_k)\\
	& & {}+2\E(w_1(1)w_1(2))^2
	\summe{i,j=1}{n}\left(f(x_i,x_j)\right)^2\nonumber
\end{eqnarray*}
and
\begin{eqnarray*}
	\E\hut{\theta_1^*}\hut{\theta_2^*} & = & 
	\hspace{0.5cm}\E\summe{i,j,k,l=1l}{n}f(x_i,x_j)f(x_k,x_l)
	w_1(i)w_1(j)w_2(k)w_2(l)\nonumber\\
	& & {}+4\E\summe{i,j,k=1}{n}f(x_i,x_j)f(x_i,x_k)
	w_1(i)w_1(j)w_2(i)w_2(k)\\
	& & {}+2\E\summe{i,j=1}{n}\left(f(x_i,x_j)\right)^2
	w_1(i)w_1(j)w_2(i)w_2(j)\\
	& = & \hspace{0.67cm}\left(\E w_1(1)w_1(2)\right)^2
	\summe{i,j,k,l=1}{n}
	f(x_i,x_j)f(x_k,x_l)\nonumber\\
	& & {}+4\left(\E w_1(1)w_1(2)\right)^2
	\summe{i,j,k=1}{n}f(x_i,x_j)f(x_i,x_k)\\
	& & {}+2\left(\E w_1(1)w_1(2)\right)^2
	\summe{i,j=1}{n}\left(f(x_i,x_j)\right)^2.\nonumber
\end{eqnarray*}
This yields
\begin{eqnarray*}
	\lim_{N\to\infty}\hut{v_N^*} & = & \hspace{0.7cm}
	\left[\E w_1(1)w_1(2)w_1(3)w_1(4)-
	\left(\E w_1(1)w_1(2)\right)^2\right]\cdot{}
	\nonumber\\
	& & \hspace{1cm}{}\cdot\summe{w,x,y,z\in\Phi}{}
	f(w,x)f(y,z)\nonumber\\
	& & {}+4\left[\E(w_1(1))^2 w_1(2)w_1(3)-
	\left(\E w_1(1)w_1(2)\right)^2\right]\cdot{}\\
	& & \hspace{1cm}{}\cdot\summe{x,y,z\in\Phi}{}f(x,y)f(x,z)\\
	& & {}+2\left[\E(w_1(1)w_1(2))^2-
	\left(\E w_1(1)w_1(2)\right)^2\right]\cdot{}\nonumber\\
	& & \hspace{1cm}{}\cdot\summe{x,y\in\Phi}{}\left(f(x,y)\right)^2\\
	& = & \alpha_4\summe{i,j,k,l=1}{n}
	f(x_i,x_j)f(x_k,x_l)\nonumber\\
	& & {}+4\alpha_3\summe{i,j,k=1}{n}
	f(x_i,x_j)f(x_i,x_k)\\
	& & {}+2\alpha_2\summe{i,j=1}{n}
	\left(f(x_i,x_j)\right)^2\nonumber
\end{eqnarray*}
with
\begin{eqnarray}\label{tron}
	\alpha_4 & = & (-4n^2+10n-6)/n^3\nonumber\\
	\alpha_3 & = & (n^3-7n^2+12n-6)/n^3\\
	\alpha_2 & = & (3n^3-11n^2+14n-6)/n^3\nonumber
\end{eqnarray}

\subsection*{Equation (\ref{result})}

The expectation value of $\hut{v^*}$ is
\begin{eqnarray*}
	\E\hut{v}^* & = & 
	\alpha_4\int\varrho^{(4)}(x_1,x_2,x_3,x_4)f(x_1,x_2)f(x_3,x_4)
	\textrm{d}\, x_1\textrm{d}\, x_2\textrm{d}\, x_3\textrm{d}\, x_4\nonumber\\
	& & {}+4\alpha_3\int\varrho^{(3)}(x_1,x_2,x_3)f(x_1,x_2)f(x_1,x_3)
	\textrm{d}\, x_1\textrm{d}\, x_2\textrm{d}\, x_3\\
	& & {}+2\alpha_2\int\varrho^{(2)}(x_1,x_2)\left(f(x_1,x_2)\right)^2
	\textrm{d}\, x_1\textrm{d}\,x_2\\
	& = & s_4\alpha_4+4s_3\alpha_3+2s_2\alpha_2
\end{eqnarray*}
In the limiting case ($n\to\infty$) it is
\begin{eqnarray*}
	\E\hut{v^*} & = & 4s_3+6s_2
\end{eqnarray*}
(see Equation (\ref{tron})).


\end{document}